\documentclass[11pt,reqno]{amsart}
\usepackage[a4paper,textwidth=153mm,textheight=235mm,centering]{geometry}
\usepackage{amsmath,amssymb,amsthm,iftex}
\usepackage[unicode,hidelinks,bookmarksnumbered=true]{hyperref}
\hypersetup{pdftitle={Rank-two twists of the Fermat cubic with 3-primary Tate-Shafarevich group (Z/3Z)\string^(2m)}}
\newcommand{\Q}{\mathbb Q}
\newcommand{\F}{\mathbb F}
\newcommand{\Z}{\mathbb Z}
\newcommand{\Gal}{\operatorname{Gal}}
\newcommand{\Sel}{\operatorname{Sel}}
\DeclareFontFamily{U}{shacyr}{}
\DeclareFontShape{U}{shacyr}{m}{n}{<-> wncyr10}{}
\DeclareSymbolFont{shafont}{U}{shacyr}{m}{n}
\DeclareMathSymbol{\shaSymbol}{\mathord}{shafont}{"58}

\ifXeTeX
  \usepackage[no-math]{fontspec}
  \ExplSyntaxOn
  \clist_map_inline:nn
    {CMU~Serif,Times~New~Roman,DejaVu~Serif,Noto~Serif,Liberation~Serif}
    {
      \IfFontExistsTF{#1}{
        \newfontface\shaCyrillicFont{#1}
        
        \clist_map_break:
      }{}
    }
  \ExplSyntaxOff
\else
  \usepackage{accsupp}
  
\fi
\DeclareRobustCommand{\Sha}{\ensuremath{\mathord{\shaGlyph}}}
\newtheorem{theorem}{Theorem}[section]
\newtheorem{proposition}[theorem]{Proposition}
\newtheorem{lemma}[theorem]{Lemma}
\newtheorem{corollary}[theorem]{Corollary}
\newtheorem{mainresult}{Theorem}

\theoremstyle{definition}
\newtheorem{definition}[theorem]{Definition}
\newtheorem{example}[theorem]{Example}
\numberwithin{equation}{section}
\title[Rank-two twists of the Fermat cubic]{%
	Rank-two twists of the Fermat cubic with
	$3$-primary Tate--Shafarevich group
	$(\mathbb{Z}/3\mathbb{Z})^{2m}$}
\author{Hongfeng Wu}
\address{College of Science, North China University of Technology, Beijing, China}
\email{whfmath@gmail.com}

\author{Li Zhu}
\address{School of Mathematical Sciences, Guizhou Normal University, Guiyang, China}
\email{lizhumath@pku.edu.cn}
\date{}
\begin{document}

\begin{abstract}
For each integer $m\ge0$ we construct infinitely many positive squarefree
integers $k$ for which the Fermat cubic $E_k\colon x^3+y^3=k$ has
Mordell--Weil rank two over $\Q$ and
$\Sha(E_k/\Q)[3^\infty]\simeq(\Z/3\Z)^{2m}$.
Fixing $2m$ prime divisors of $k$ congruent to $-1$ modulo nine,
with product $D$, we obtain
$\gg_D X^{1/3}/(\log X)^6$ distinct parameters up to $X$.
Two rational points on a $3$-isogenous curve come from the Stewart--Top
sextic. Modulo their Kummer classes and a torsion class, the isogeny
Selmer group has a basis given by the fixed primes. We realize every
alternating matrix in this basis by prescribing the Frobenius classes
of six primes in a fixed Galois extension. A nondegenerate matrix
forces the rank to be exactly two and determines the entire primary group.
\end{abstract}
\maketitle

\section{Introduction}\label{sec:intro}
Let $E_k$ be the smooth plane cubic $x^3+y^3=kz^3$ with origin
$[1:-1:0]$, where $k$ is a positive cube-free integer. What finite
groups can occur as $\Sha(E_k/\Q)[3^\infty]$ when the
Mordell--Weil rank is fixed at two? Cassels proved that
$\dim_{\F_3}\Sha(E_k/\Q)[3]$ is unbounded in this cubic twist family
\cite{CasVI}. Stewart and Top \cite[Theorem~7 and Sec.~9]{ST}
constructed $\gg X^{1/3}$ cube-free twists with parameter of
absolute value at most $X$ and Mordell--Weil rank at least two.
On a subfamily of their construction we obtain rank exactly two and
prescribe the entire $3$-primary group.

\begin{mainresult}\label{thm:main}
Fix an integer $m\ge0$ and distinct primes
$q_1,\ldots,q_{2m}\equiv-1\pmod9$. Write $D=\prod_i q_i$,
with $D=1$ when $m=0$. There is a family $\mathcal K_{2m}(D)$ of
positive squarefree integers $k\equiv1\pmod9$ such that every
$k\in\mathcal K_{2m}(D)$ satisfies
\[
  \operatorname{rank}E_k(\Q)=2,\qquad
  \Sha(E_k/\Q)[3^\infty]\simeq(\Z/3\Z)^{2m} .
\]
Every $k$ is $D$ times six distinct primes prime to $6D$, two
congruent to $4$ and four congruent to $5$ modulo nine. For all
sufficiently large $X$,
\[
  \#\{k\le X:k\in\mathcal K_{2m}(D)\}
    \gg_D \frac{X^{1/3}}{(\log X)^6}.
\]
Distinct parameters give nonisomorphic curves over both $\Q$ and
$K=\Q(\sqrt{-3})$.
\end{mainresult}

The parity restriction on the dimension is necessary. On a finite $3$-primary group the
Cassels--Tate pairing is perfect and alternating, so an elementary
abelian group carrying it has even dimension. Theorem~\ref{thm:main}
therefore exhausts the finite elementary abelian possibilities for
the $3$-primary group at rank two. Both the
implied constant and the threshold for $X$ may depend on $D$.

Kezuka's lower bound for this family clarifies the choice of prime
factors. Write $\nu_{\rm in}(k)$ for the number of distinct prime
divisors of $k$ inert in $K$. At $k>1$ with $k\equiv1\pmod9$,
\cite[Theorem~1.6]{Kez} gives
\[
 \dim_{\F_3}\Sha(E_k/\Q)[3]\ \ge \nu_{\rm in}(k)-2-\operatorname{rank}E_k(\Q).
\]
Our parameters have $2m+4$ inert prime divisors. At rank two the right
side is $2m$. Thus Theorem~\ref{thm:main} attains equality in Kezuka's
bound and shows that every $3$-primary class on these curves is killed
by three.

\subsection*{Realizing the pairing on fixed classes}
Let $A_k\colon y^2=x^3+16k^2$ be the $3$-isogenous partner of $E_k$.
Write $\phi:E_k\to A_k$ for the isogeny with kernel $\mu_3$ and
$\widehat\phi:A_k\to E_k$ for its dual, with kernel generated by
$T=(0,4k)$. The group $U_k=\Sel_\phi(E_k/\Q)$ is a subspace of
$\Q^\times/\Q^{\times3}$. Denote by $\delta_k:A_k(\Q)\to U_k$
the Kummer map, and by $\Theta_k$ the pullback of the Cassels--Tate
pairing under $U_k\to\Sha(E_k/\Q)[\phi]$. Its values lie in
$\tfrac13\Z/\Z$, identified with $\F_3$ by $j/3\mapsto j$.
The classes of rational points on $A_k$ lie in its radical.

\begin{mainresult}\label{thm:realization}
Fix $m,q_1,\ldots,q_{2m},D$ as in Theorem~\ref{thm:main}, and let $B$
be any alternating $2m\times 2m$ matrix over $\F_3$. For all sufficiently
large $X$ there are $\gg_D X^{1/3}/(\log X)^6$ distinct positive
squarefree $k\le X$ of the prime factorization specified there with
the following properties. There exist $P,Q\in A_k(\Q)$ independent
modulo torsion. The space
\[
 W_k=\langle\delta_k(T),\delta_k(P),\delta_k(Q)\rangle
\]
has dimension three, and $\Sel_{\widehat\phi}(A_k/\Q)=0$.
In the ordered basis of $U_k/W_k$ given by the cube classes
$x_i=[q_i]$, the induced pairing has matrix $B$.
\end{mainresult}

The ordering matters. The entry $B_{ij}$ specifies the pairing between
the classes of the two particular primes $q_i$ and $q_j$. Singular
matrices occur as well. In Section~\ref{sec:criterion} we show that,
writing $d=\operatorname{rank} B$, the number of cyclic factors of
exact order three in the quotient of $\Sha(E_k/\Q)[3^\infty]$ by
its maximal divisible subgroup equals $d$. The rational rank lies
between $2$ and $2+2m-d$. When $d=2m$, these bounds force rank two.
The pairing is then nondegenerate on $\Sha(E_k/\Q)[3]$. If $y$ had
order nine, its nonzero triple would pair trivially with every
three-torsion class by $\langle3y,x\rangle=\langle y,3x\rangle$.
This explains how the pairing determines the exponent without an
assumption of finiteness.

\subsection*{Construction and relation to previous work}
Two identities among the six linear factors of the Stewart--Top sextic
give the points $P,Q$. We keep the factors linear in two parameters,
and impose cubic residue conditions under which the descent matrix
has three independent rows and zero columns at the $q_i$. This leaves
the classes $[q_i]$ free after the rational Kummer classes have been
accounted for. The descent uses the local matrix of Elkies and Rogers
\cite[Sec.~3]{ER} and the isogeny pairing formalism of van Beek
and Fisher \cite{vBF}.

Koymans and Smith \cite[Definition~1.4 and Theorem~1.5]{KS} determine
the distribution of a modified $3$-Selmer rank in cubic twist families. Their Proposition~4.1 compares pairings
when one prime is replaced. Here the six prime factors are values of
linear forms in the same two integers. We compare their product twist
directly with the twist by $D$. For each pair of fixed classes, the
change in its pairing entry is a sum of local invariants. These sums
are read from central Frobenius coordinates in a single extension
chosen for $D$. Assigning all the coordinates of one central element
prescribes all the entries at once. The prime-values theorem of
Frei and Tinkov\'a \cite[Theorem~1.4 and Corollary~1.6]{FT} realizes the
six Frobenius conditions simultaneously and gives the count in
Theorem~\ref{thm:realization}.

For comparison, in quadratic twist families Bhargava, Klagsbrun,
Lemke Oliver and Shnidman \cite{BKLOS} use independent isogenies
to force elements of given prime order in Tate--Shafarevich groups
for positive proportions of twists in suitable families. Exact
rank-two constructions include Zywina's varying-$j$ family
\cite{Zywina} and the quadratic twists of Achenjang and Woo \cite{AW}.
For the Fermat cubic, Byeon and Jeong
\cite[Theorem~1.1 and Proposition~3.1]{BJ} study parameters $pq$
with $p,q\equiv-1\pmod9$ and obtain exact rank two on their family
under the Mordell--Weil parity conjecture.

Sections~\ref{sec:descent}--\ref{sec:both} compute the isogeny
Selmer groups and the rational classes. Sections~\ref{sec:group}
and~\ref{sec:counting} prove Theorems~\ref{thm:realization}
and~\ref{thm:main}. Example~\ref{ex:nonzero-fiber} gives a single
fiber with $\Sha(E_k/\Q)[3^\infty]\simeq(\Z/3\Z)^2$.
Section~\ref{sec:six-split} gives a disjoint family on which the
rational classes alone imply vanishing. In Section~\ref{sec:consequences},
base change to the CM field doubles the Mordell--Weil rank, and
restriction identifies the $3$-primary Tate--Shafarevich groups over $\Q$ and $K$.

\section{Preliminaries}\label{sec:prelim}

\subsection{Models and isogenies}
Let $k$ be a positive cube-free integer. The plane model
$E_k\colon a^3+b^3=kc^3$ has origin $[1:-1:0]$. On $c=1$ the map
\[
  (a,b)\longmapsto
  \left(\frac{12k}{a+b},\,\frac{36k(b-a)}{a+b}\right)
\]
gives the short Weierstrass equation $Y^2=X^3-432k^2$.
Write $A_k\colon V^2=U^3+16k^2$. The point $T=(0,4k)$ has order
three, and quotienting $A_k$ by $\langle T\rangle$ gives $E_k$
\cite{CP,ER}. Our notation for the two degree-three isogenies is
\[
  E_k\xrightarrow{\ \phi\ } A_k
       \xrightarrow{\ \widehat\phi\ } E_k,\qquad
  \widehat\phi\circ\phi=[3].
\]
In particular $A_k[\widehat\phi]\simeq\Z/3\Z$ and the Weil pairing
identifies $E_k[\phi]$ with $\mu_3$, the group of cube roots of unity.
If $3\nmid k$, the displayed equation for $E_k$ has discriminant of
$3$-adic valuation $9$. It is minimal at $3$, and its reduction there
is cuspidal.

For $k\ne0$ the model has $j=0$. Over $K=\Q(\zeta_3)$ the map
$(X,Y)\mapsto(\zeta_3X,Y)$ realizes complex multiplication by
$\mathcal O_K=\Z[\zeta_3]$. All parameters in Theorem~\ref{thm:main}
satisfy $k\equiv1\pmod9$. On $E_k$ the Kodaira symbol at three
is then $\mathrm{IV}^*$ and the Tamagawa number is three
\cite[Lemma~2.1, Table~1]{FujN}.

\subsection{Selmer groups}
Let $L$ be a number field, let $L_v$ be its completion at $v$, and
write $H^1(L,-)$ for continuous Galois cohomology. For an elliptic
curve $B/L$ put
\[
 \Sha(B/L)=
 \ker\!\left(H^1(L,B)\longrightarrow\prod_v H^1(L_v,B)\right).
\]
For an isogeny $u\colon B\to B'$ over $L$, the connecting map
$\delta_{u,v}$ sends $B'(L_v)/uB(L_v)$ into $H^1(L_v,B[u])$.
The $u$-Selmer group $\Sel_u(B/L)$ is the subgroup of $H^1(L,B[u])$
whose restriction at every $v$ lies in the image of $\delta_{u,v}$.
The global and local Kummer sequences give \cite[Chap.~X]{Sil}
\begin{equation}\label{eq:kummer-pre}
 0\longrightarrow B'(L)/uB(L)\longrightarrow \Sel_u(B/L)
 \longrightarrow\Sha(B/L)[u]\longrightarrow0.
\end{equation}
Here $\Sha(B/L)[u]$ is the kernel of $u$ on $\Sha(B/L)$. For $u=[3]$
we write $\Sel_3(B/L)$. In particular,
\begin{equation}\label{eq:dimensions}
 \dim_{\F_3}\Sel_3(B/L)=\operatorname{rank}B(L)
 +\dim_{\F_3}B(L)[3]+\dim_{\F_3}\Sha(B/L)[3].
\end{equation}
For a pair of dual isogenies, local Tate duality identifies the local
condition for one with the annihilator of the other under the local
Tate pairing \cite[Sec.~3]{ER}. We use this to compute both Selmer
groups from one matrix.

Throughout the paper set $U_k=\Sel_\phi(E_k/\Q)$,
$V_k=\Sel_{\widehat\phi}(A_k/\Q)$,
$G_k=\Sha(E_k/\Q)[3^\infty]$ and $G'_k=\Sha(A_k/\Q)[3^\infty]$.
We omit $k$ from $U_k,V_k$ when the parameter is fixed.
For an abelian group $G$, write $G[3^\infty]=\bigcup_{a\ge1}G[3^a]$.
Vanishing of $G[3]$ implies vanishing of $G[3^\infty]$.

\subsection{Cubic characters and Frobenius conditions}
Write $[a]$ for a cube class in $\Q^\times/\Q^{\times3}$ and $v_\ell$
for the valuation normalized by $v_\ell(\ell)=1$. At a prime
$p\equiv1\pmod3$ we use cubic characters in additive notation, as
surjections $\F_p^\times\to\F_3$. For a fixed primitive cube root
$\omega\in\F_p$, the normalization is $a^{(p-1)/3}=\omega^{\chi(a)}$.
Changing $\omega$ to $\omega^{-1}$ negates $\chi$. We use arithmetic
Frobenius, which raises residues to the $p$th power at an unramified
rational prime $p$.

We recall the form of the prime-values theorem needed below
\cite[Theorem~1.4 and Corollary~1.6]{FT}. Let $H/\Q$ be a fixed
finite Galois extension and $C_1,\ldots,C_r$ conjugacy classes of its
Galois group. Let $f$ be the least positive integer such that the
maximal abelian subextension $H^{\rm ab}$ is contained in
$\Q(\zeta_f)$. Let $\psi_1,\ldots,\psi_r$ be fixed integral affine
forms in two variables with nonzero, pairwise nonproportional linear
parts. Positivity of the finite factors in the theorem follows from
the following two conditions.
\begin{enumerate}
\item For each prime $\ell$ some residue pair modulo $\ell$ makes all
$\psi_i$ nonzero modulo $\ell$.
\item Some common residue pair modulo $f$ makes each
$b_i=\psi_i(x,y)$ a unit modulo $f$, with
$\zeta_f\mapsto\zeta_f^{b_i}$ agreeing on $H^{\rm ab}$ with $C_i$.
\end{enumerate}
If $\Omega\subset(-1,1)^2$ is a fixed closed rectangle of positive
area on which all the linear parts are positive, the weighted count
on $B\Omega$ then has a positive main term of order $B^2$.
Sections~\ref{sec:counting} and~\ref{sec:six-split} verify these conditions for the two extensions used below.
For six forms we will use the following consequence on a dilated rectangle
$Y\Omega$. When the weighted count is $\gg Y^2$ and the six values
are each $\asymp Y$, tuples of six distinct primes contribute
$\gg Y^2/(\log Y)^6$ parameter pairs. Indeed each form assumes
a prescribed value at $O(Y)$ lattice points and there are
$O(Y^{1/2}\log Y)$ proper prime powers in its range. The total
weight of tuples involving such a power is
$O(Y^{3/2}(\log Y)^7)$. Coincident values occur on finitely many lines
of $O(Y)$ lattice points. Each remaining tuple has weight
$O((\log Y)^6)$. All constants may depend on the fixed forms and fields.

\section{Cubic isogeny descent on the Fermat cubic}\label{sec:descent}

For the descent we use the following plane model of the isogeny partner
$A_k$ from Section~\ref{sec:prelim}.
\[
  uv(u+v)=k,\qquad V^2=U^3+16k^2,
  \qquad U=4k/v,\quad V=4k(2u+v)/v .
\]
The zero on the plane model $uv(u+v)=kw^3$ is $[1:0:0]$.
The tangent function at $T$ gives a connecting map with values in
$H^1(\Q,\mu_3)=\Q^\times/\Q^{\times3}$. With our choice of $T$, it is
\begin{equation}\label{eq:cubic-kummer}
  \delta(R)=[V(R)-4k]\quad(R\ne O,T),\qquad
  \delta(O)=[1],\qquad \delta(T)=[k]^{-1}.
\end{equation}
Indeed, the divisor of $V-4k$ is $3(T)-3(O)$, and
$(V-4k)(V+4k)=U^3$ gives the value $[(8k)^{-1}]=[k]^{-1}$
at $T$ \cite{CP}. In particular, $\delta(-T)=[k]$.
For a finite point on the $uv$ model the formula becomes
$\delta(u,v)=[ku/v]$.

\begin{proposition}\label{prop:descent}
Let $k>1$ be odd and squarefree with $k\equiv1\pmod 9$.
List its prime factors as $\ell_1,\ldots,\ell_t$, and write
$p_1,\ldots,p_h$ for those congruent to $1\pmod 3$.
Choose a nonzero additive cubic character
$\chi_i\colon\F_{p_i}^{\times}\longrightarrow\F_3$ at each $p_i$.
Let $\chi_9\colon(\Z/9\Z)^\times\longrightarrow\F_3$ be the character
with $\chi_9(4)=1$. Define the $(h+1)\times t$ matrix $M$ over $\F_3$ by
\begin{align*}
 M_{ij}&=
 \begin{cases}
   \chi_i(\ell_j),&\ell_j\ne p_i,\\
   -\chi_i(k/p_i),&\ell_j=p_i
 \end{cases}
 &&(1\le i\le h),&
 M_{h+1,j}&=\chi_9(\ell_j).
\end{align*}
Then $\Sel_\phi(E_k/\Q)\simeq\ker M$ and
$\Sel_{\widehat{\phi}}(A_k/\Q)\simeq\ker M^{\mathsf T}$.
Their dimensions are $t-\operatorname{rank}M$ and
$h+1-\operatorname{rank}M$, respectively. If $M$ has full row rank, then
\[
  \Sel_3(E_k/\Q)\simeq
  \ker M/\langle(1,\ldots,1)\rangle .
\]
\end{proposition}

\begin{proof}
We recall the local calculation behind the descent matrix of
Elkies and Rogers \cite[Sec.~3]{ER}. A class $[d]$ belongs to the
$\phi$-Selmer group precisely when the diagonal cubic
\[
  d^{-1}x^3+d y^3=kz^3
\]
has a point over every completion of $\Q$. This equation follows by
writing $V-4k$ and $V+4k$ as $d$ and $d^{-1}$ times cubes, then
taking their difference $8k$ and absorbing $8$ and the sign into the
cube variables. At a prime $\ell\nmid k$,
the three terms have distinct valuations modulo $3$ if
$3\nmid v_\ell(d)$, which rules out a point. At $\ell\nmid3k$
this is the only condition. After removing cubes, a unit $d$ gives
a smooth plane cubic modulo $\ell$. The genus-one Weil bound gives a
point and Hensel's lemma lifts it. The real place gives no condition.
Hence each candidate is represented uniquely modulo cubes by
$d=\prod_j\ell_j^{e_j}$, with exponents taken in $\{0,1,2\}$.

At $\ell\mid k$, the local image contains $[k]=\delta(-T)$.
Dividing by a power of this class reduces the question to a unit $d$.
In primitive integral coordinates, reduction of
$x^3+d^2 y^3=kd z^3$ forces $x$ and $y$ to be units.
Indeed, if both vanished modulo $\ell$, the two sides would have
incompatible valuations. Thus $d$ is a cube modulo $\ell$.
Conversely such a unit is a local cube, so gives a point with $z=0$.
The local image is therefore $\langle[k]\rangle$.
At $p_i=\ell_j$ the unit test on $d/k^{e_j}$ gives the $i$th row of
$M$. At primes congruent to $2\pmod3$ every unit is a cube.

At $3$ a unit $d$ is a cube exactly when $d\equiv\pm1\pmod9$.
If $d$ is a noncube unit, the coefficients $d^{-1},d,-k$ are,
up to sign and order, $1,2,4$ modulo $9$. In primitive coordinates
each cube is $0$ or $\pm1$ modulo $9$, with at least one unit.
No nonempty signed sum of distinct elements of $\{1,2,4\}$ is zero
modulo $9$. Thus the local image at $3$ is trivial. This gives the
last row of $M$ and proves the right-kernel formula.

For the dual isogeny the global classes are cubic characters.
By local Tate duality, a dual class annihilates the local cubic
Kummer image under the local reciprocity pairing. At $\ell\nmid 3k$
this image consists of the unit classes. Thus a character must be unramified away from
$3k$ and vanish on $k$ at each $\ell\mid k$. There is no condition
at $3$. By Kronecker--Weber, the characters unramified away from
$3k$ have basis $\chi_1,\ldots,\chi_h,\chi_9$, viewed as Dirichlet
characters. With arithmetic Frobenius for uniformizers, the local
value of $\chi_i$ on $k$ at $\ell_j\ne p_i$ is $\chi_i(\ell_j)$.
At $p_i$ it is $-\chi_i(k/p_i)$, since local units act by inverse
powers on the $p_i$th roots of unity. The character $\chi_9$ takes
the value $\chi_9(\ell_j)$ on $k$ at $\ell_j$. Thus the dual local
conditions give precisely the left kernel of $M$.

Finally, $E_k(\Q)[3]=0$. The nonzero geometric $3$-torsion points
have either $X=0$, with $Y^2=-432k^2$, or $X^3=12^3k^2$, which is
impossible over $\Q$ by the valuation at any prime of $k$.
The factorization $[3]=\widehat{\phi}\circ\phi$ gives the exact sequence
\begin{equation}\label{eq:sel-sequence}
0\longrightarrow \langle T\rangle
 \xrightarrow{\delta} \Sel_\phi(E_k/\Q)
 \longrightarrow \Sel_3(E_k/\Q)
 \longrightarrow \Sel_{\widehat{\phi}}(A_k/\Q).
\end{equation}
The first image is generated by $\delta(-T)=[k]$, whose exponent
vector is $(1,\ldots,1)$. Full row rank gives a trivial last group
and proves the quotient formula.
\end{proof}

\section{A criterion for the entire three-primary group}\label{sec:criterion}

Throughout this section assume $k>1$ is cube-free. The division polynomial
calculation in the proof of Proposition~\ref{prop:descent} gives
$E_k(\Q)[3]=0$ for all such $k$.
Write \(f=\phi:E_k\to A_k\) and \(g=\widehat\phi:A_k\to E_k\).
We use the same letters for their maps on
\(G_k=\Sha(E_k/\Q)[3^\infty]\) and \(G'_k=\Sha(A_k/\Q)[3^\infty]\).
The isogeny Kummer sequence gives a surjection
\(U\to G_k[f]\) with kernel \(A_k(\Q)/fE_k(\Q)\).
Let \(\Theta\) be the pullback along this surjection of the
Cassels--Tate pairing on \(G_k\). Its values lie in
\(\tfrac13\Z/\Z\), which we identify with \(\F_3\) by \(j/3\mapsto j\).
In particular the classes of rational points on \(A_k\) are in its radical.

There is also a descent interpretation of the radical. The coefficient
sequence \(0\to A_k[g]\to A_k[3]\xrightarrow{g}E_k[f]\to0\)
gives the sequence dual to \eqref{eq:sel-sequence}. Since
\(E_k[f](\Q)=\mu_3(\Q)=0\), its first map on Selmer groups is
injective. Its last image is the radical of \(\Theta\)
\cite[Introduction, pairing (2), and Definition~1.2]{vBF}. Thus
\begin{equation}\label{eq:partner-radical}
  0\longrightarrow V\longrightarrow\Sel_3(A_k/\Q)
    \xrightarrow{g_*}\operatorname{rad}\Theta\longrightarrow0 .
\end{equation}

\begin{lemma}\label{lem:general-rank-bound}
Put \(d=\operatorname{rank}_{\F_3}\Theta\) and
\(r=\operatorname{rank}E_k(\Q)\). For any \(V\),
\[
 d\leq\dim_{\F_3}G_k[3]=\dim_{\F_3}\Sel_3(E_k/\Q)-r .
\]
If \(r\geq2\) and \(d=\dim\Sel_3(E_k/\Q)-2\), then \(r=2\) and
\(G_k\simeq(\Z/3\Z)^d\).
\end{lemma}

\begin{proof}
The rank of the pairing on the subgroup \(G_k[f]\) is at most the
dimension of \(G_k[3]\). Kummer theory gives the equality since
\(E_k(\Q)[3]=0\). Under the extra hypotheses the inequalities force
\(r=2\) and \(G_k[f]=G_k[3]\), with nondegenerate restricted pairing.
Adjointness of multiplication by \(3\) shows that
\(G_k[3]\cap 3G_k\) lies in the radical. An element of order nine
would have a nonzero triple in this intersection. Thus every element
of \(G_k\) is killed by three.
\end{proof}

For a cofinitely generated \(3\)-primary group \(G\), let \(\mu_1(G)\)
count the cyclic direct summands of order exactly three in the finite
quotient by its maximal divisible subgroup.

\begin{theorem}\label{thm:criterion}
Suppose \(V=0\) and \(E_k(\Q)\) has two independent points. Write
\(a=\dim_{\F_3} U\) and \(d=\operatorname{rank}_{\F_3}\Theta\).
Then \(g\) identifies \(G'_k\) with \(3G_k\), and
\[
  \mu_1(G_k)=d,\qquad 2\leq r\leq a-1-d .
\]
Moreover,
\[
 d=a-3 \quad\Longleftrightarrow\quad
 (r=2\ \text{and}\ G'_k=0) \quad\Longleftrightarrow\quad
 (r=2\ \text{and}\ G_k\simeq(\Z/3\Z)^{a-3}).
\]
\end{theorem}

\begin{proof}
The isogeny Kummer sequence for \(g\) and the assumption \(V=0\)
give \(G'_k[g]=0\). To prove that \(f:G_k\to G'_k\) is surjective,
write \(D,D'\) for the maximal divisible subgroups of \(G_k,G'_k\).
These groups map onto each other under \(f,g\), since both composites
are multiplication by \(3\), which is surjective on a divisible group.
The induced map \(\bar g:G'_k/D'\to G_k/D\) is injective. Indeed,
if \(g(x)\in D\), choose \(y\in D'\) with \(g(y)=g(x)\), and use
\(G'_k[g]=0\). The two quotients are finite and carry perfect
Cassels--Tate pairings for which \(\bar f\) and \(\bar g\) are
adjoint. Thus \(\bar f\) is surjective. Together with \(f(D)=D'\),
this gives \(f(G_k)=G'_k\), and hence \(g(G'_k)=gf(G_k)=3G_k\).

By \eqref{eq:sel-sequence}, \(\dim\Sel_3(E_k/\Q)=a-1\). By
\eqref{eq:partner-radical}, \(\dim\Sel_3(A_k/\Q)=a-d\).
The two curves have the same rank, and their rational \(3\)-torsion
has dimension zero and one, respectively. Their Kummer sequences give
\[
 \dim G_k[3]=a-1-r,\qquad
 \dim G'_k[3]=a-1-r-d .
\]
Since \(G'_k\simeq3G_k\), subtracting these dimensions counts the
cyclic factors of order exactly three. Every longer cyclic factor
and every copy of \(\Q_3/\Z_3\) contributes one to each dimension.
The second dimension gives the rank bound.

If \(d=a-3\), the bound gives \(r=2\), and the second dimension gives
\(G'_k[3]=0\), hence \(G'_k=0\). At rank two, \(G'_k=0\) is
equivalent to \(3G_k=0\), in which case the first dimension gives
\(G_k\simeq(\Z/3\Z)^{a-3}\). Conversely, this elementary abelian group
has \(\mu_1(G_k)=a-3\).
\end{proof}

For the odd squarefree parameters \(k\equiv1\pmod9\) in
Proposition~\ref{prop:descent}, let \(t\) count the prime factors of
\(k\) and \(h\) count those splitting in \(\Q(\sqrt{-3})\). If the
descent matrix has full row rank \(h+1\), then that proposition gives \(V=0\) and \(a=t-h-1\). In particular, in
Theorem~\ref{thm:criterion} the number \(a-3=t-h-4\) is the number
of inert prime factors of \(k\) minus four.
Appendix~\ref{app:1729} illustrates the rank bound when the dual isogeny Selmer group is nonzero.

\section{Rational points and fixed Selmer classes on the Stewart--Top sextic}\label{sec:both}

\subsection{Two independent rational points}
Write the sextic of Stewart and Top \cite[Sec.~9]{ST} as twice the product
of the linear forms
\[
 (L_1,L_2,L_3,L_4,L_5,L_6)=(n_1,n_2,n_1-n_2,n_1+n_2,2n_1+n_2,n_1+2n_2).
\]
In both constructions below, $L_i=d_i p_i$ with signed integer $d_i$
and distinct positive primes $p_i$ coprime to every $d_j$. Write $D$ for
the product of the fixed primes, which are distinct from the $p_i$, and
put $k=D\prod_i p_i$. The signed factorizations satisfy
\begin{equation}\label{eq:sextic-points}
 2\prod_{i=1}^6L_i=\varepsilon\,60^3 k,\qquad \varepsilon\in\{1,-1\}.
\end{equation}

\begin{proposition}\label{prop:points}
Under these hypotheses $A_k(\Q)$ contains two points $P,Q$ whose
images under $\widehat\phi$ are independent modulo torsion on $E_k$.
The valuations modulo three on the moving primes are
\begin{equation}\label{eq:point-vectors}
\begin{array}{c|rrrrrr}
 &p_1&p_2&p_3&p_4&p_5&p_6\\ \hline
 [k]           &1&1&1&1&1&1\\
 \delta(P)/[k] &1&0&1&2&0&2\\
 \delta(Q)/[k] &0&2&2&1&1&0
\end{array}
\end{equation}
\end{proposition}

\begin{proof}
On the model $uv(u+v)=k$ of $A_k$, take the coordinates
\[
\begin{array}{c|rrr}
 &60\varepsilon u&60\varepsilon v&60\varepsilon(u+v)\\ \hline
 P&L_1L_3&-L_6L_4&-2L_2L_5\\
 Q&L_5L_4& L_2L_3& 2L_1L_6 .
\end{array}
\]
Expansion in $n_1,n_2$ gives $L_1L_3-L_6L_4=-2L_2L_5$ and
$L_5L_4+L_2L_3=2L_1L_6$. The three entries in each row have
product $2\prod_i L_i$, giving $uv(u+v)=k$
by \eqref{eq:sextic-points}. The formula
$\delta(u,v)=[ku/v]$ from \eqref{eq:cubic-kummer} yields
\[
 \delta(-T)=[k],\qquad
 \delta(P)/[k]=[-L_1L_3/(L_6L_4)],\qquad
 \delta(Q)/[k]=[L_5L_4/(L_2L_3)].
\]
Taking valuations gives \eqref{eq:point-vectors}. Its minor on
$p_1,p_2,p_5$ has determinant $1$ in $\F_3$.
The Weil pairing gives $A_k(\Q)[3]=\langle T\rangle$, so the
rational torsion subgroup modulo three has dimension one.
Its image under $\delta$ is generated by the nontrivial class $[k]$.
A primitive integral relation between $P$ and $Q$ modulo torsion would
therefore give a nontrivial relation among the three rows modulo three.
The finite kernel of $\widehat\phi$ preserves independence modulo torsion.
\end{proof}

Let $W_k\subseteq\Sel_\phi(E_k/\Q)$ be the three-dimensional span of
$[k],\delta(P)/[k],\delta(Q)/[k]$. Each class is rational and hence
lies in the radical of the isogeny pairing.

In both patterns $k$ is odd and squarefree with $k\equiv1\pmod9$.
Order the columns of the matrix $M$ of Proposition~\ref{prop:descent}
by the fixed primes $q_j$ followed by $p_1,\ldots,p_6$. Every $q_j$ is
$-1\pmod9$, and every $p_i$ is $4$ or $-4\pmod9$. Hence the row at $3$
is $w=(0,\ldots,0\mid1,1,1,1,1,1)$.

\begin{lemma}\label{lem:eliminate}
Write $L_j=a_j n_1+b_j n_2$ and $v_j=(a_j,b_j)$. For each split $p_i$
the corresponding row of $M$ can be replaced by a row whose entry
at $q_j$ is $\chi_i(q_j)$ and whose entry at $p_j$ for $j\ne i$ is
\[
 \chi_i\bigl(\det(v_i,v_j)\bigr)-\chi_i(d_j).
\]
Its entry at $p_i$ is the negative of the sum of its other entries.
\end{lemma}

\begin{proof}
The distinctness of the $p_j$ and the coprimality to the $d_j$ give
$(n_1,n_2)\ne(0,0)$ modulo $p_i$. Since $L_i=0$ there,
$(n_1,n_2)=h_i(-b_i,a_i)$ in $\F_{p_i}^2$ for some $h_i\ne0$.
Thus $p_j\equiv h_i\det(v_i,v_j)/d_j\pmod{p_i}$ for $j\ne i$.
Subtracting $\chi_i(h_i)w$ removes the common term from these entries.
Both the original row and $w$ have sum zero.
\end{proof}

The same three rows in the moving columns will occur in both matrices.
Set
\begin{equation}\label{eq:moving-block}
 C=\begin{pmatrix}
 2&1&0&0&1&2\\
 0&0&1&0&1&1\\
 1&1&1&1&1&1
 \end{pmatrix}.
\end{equation}
The minor on $p_2,p_3,p_4$ is lower triangular with diagonal $(1,1,1)$.

\subsection{Two split moving primes and fixed Selmer classes}\label{sec:family}

Fix an integer $m\ge0$ and distinct primes
$q_1,\ldots,q_{2m}\equiv-1\pmod9$. Set $D=\prod_{j=1}^{2m} q_j$
and $c_0=54000D=2^4 3^3 5^3D$, with $D=1$ for $m=0$. Take
positive integers $a\equiv5$, $b\equiv4\pmod9$ with $b<c_0a/2$.
Require that
\begin{equation}\label{eq:fam-primes}
 (p_1,p_2,p_3,p_4,p_5,p_6)
 =(a,b,c_0a+b,c_0a-b,2c_0a-b,(c_0/2)a-b)
\end{equation}
are distinct primes coprime to $30D$. Only $p_2,p_3$ split in
$\Q(\zeta_3)$, since they are $4\pmod9$ and the others are
$-4\pmod9$. At the two split primes require nonzero additive cubic
characters with
\begin{equation}\label{eq:fam-chars}
\begin{array}{c|cc|c}
 &2&3&q_j\ (1\le j\le 2m)\\ \hline
 \chi_2&1&0&0\\
 \chi_3&0&1&0
\end{array}
\end{equation}
and put $t=\prod_i p_i$, $k=Dt$. Since $D\equiv(-1)^{2m}=1\pmod9$ and $t\equiv4^2(-4)^4\equiv1\pmod9$, we have $k\equiv1\pmod9$.
With $n_1=-c_0a$, $n_2=b$, the signed factors are
$(d_1,\ldots,d_6)=(-c_0,1,-1,-1,-1,-2)$ and
\eqref{eq:sextic-points} holds with $\varepsilon=-1$.

\begin{proposition}\label{prop:first}
For every parameter satisfying \eqref{eq:fam-primes} and
\eqref{eq:fam-chars},
\[
 \dim_{\F_3}\Sel_\phi(E_k/\Q)=2m+3,\qquad
 \Sel_{\widehat\phi}(A_k/\Q)=0,\qquad
 \dim_{\F_3}\Sel_3(E_k/\Q)=2m+2 .
\]
For the points of Proposition~\ref{prop:points}, a basis of
$\Sel_\phi(E_k/\Q)$ is
\[
 [k],\quad \delta(P)/[k],\quad \delta(Q)/[k],\quad
 x_1=[q_1],\ldots,x_{2m}=[q_{2m}].
\]
When $m=0$ the common rational rank of the curves is two and
$\Sha(E_k/\Q)[3^\infty]=0$.
\end{proposition}

\begin{proof}
At $p_2$ the determinants for $j=1,3,4,5,6$ are
$-1,-1,-1,-2,-1$. The valuations of the $d_j$ at $2$ give
$\chi_2(d_1)=1$, $\chi_2(d_6)=1$ and zero for the other three.
Lemma~\ref{lem:eliminate} gives the off-diagonal entries $2,0,0,1,2$.
At $p_3$ the determinants for $j=1,2,4,5,6$ are $1,1,2,3,3$.
All five $d_j$ have trivial character there, since the $3$-adic and
$5$-adic valuations of $c_0$ are divisible by three and $\chi_3$
vanishes on $2$ and the $q_j$. This gives $0,0,0,1,1$ off the
diagonal. Both diagonal entries are $1$, and \eqref{eq:fam-chars}
makes the fixed-prime block zero. Hence
\begin{equation}\label{eq:fam-matrix}
 M\ \sim\ \bigl[\,0_{3\times 2m}\mid C\,\bigr].
\end{equation}
The minor of \eqref{eq:moving-block} shows that $M$ has rank three.
The $x_j$ are the unit vectors in the fixed columns. In the moving
columns the three rational classes are independent by
\eqref{eq:point-vectors}. We have therefore exhibited $2m+3$
independent classes in the right kernel, which is its full dimension.
The images of $x_1,\ldots,x_{2m}$ give a basis modulo $W_k$.
The left kernel is zero. Proposition~\ref{prop:descent} identifies
the full three-Selmer group on $E_k$ with the right kernel modulo
$\langle[k]\rangle$. At $m=0$ Proposition~\ref{prop:points} accounts
for both remaining dimensions by rational points, giving the last
assertion.
\end{proof}

\section{Fixed fields and variation of the Cassels--Tate pairing}\label{sec:group}

\subsection{The governing extension}
Retain the fixed primes $q_i\equiv-1\pmod 9$, their product $D$, and
the Selmer classes $x_i=[q_i]$. Put $K=\Q(\zeta)$ and $L=\mu_3$,
where $\zeta$ is a primitive cube root of unity. Use the additive
coordinate $\zeta\mapsto1$ on $L$ and the resulting identification
$L\otimes L\simeq\F_3$. Let $f_i$ represent $x_i$ by the Kummer
cocycle of the positive real cube root of $q_i$. For complex
conjugation $c$ we have $f_i(c)=0$. Put $W=\F_3^{2m}$ and
$Z=\F_3^{\binom{2m}{2}}$, with coordinates of $Z$ indexed by $i<j$.
Complex conjugation acts by $-1$ on $W$ and by $+1$ on $Z$.

\begin{lemma}\label{lem:group}
There is a finite Galois extension $H/\Q$, unramified away from $3D$ at
finite primes, with group $P\rtimes\langle c\rangle$, where as a set
$P=W\times Z$, and
\[
 (v,z)(v',z')=(v+v',\,z+z'+(v_i v'_j)_{i<j}),\qquad
 c(v,z)c^{-1}=(-v,z).
\]
The projection of the Galois action to $W$ over $K$ has Kummer
coordinates $f_i$. There are continuous $1$-cochains $h_{ij}$ on
$\Gal(H/\Q)$ with values in the trivial module $\F_3$ satisfying
$dh_{ij}=f_i\smile f_j$. On any subgroup where the $f_i$ vanish as
cocycles, $h_{ij}$ is a character. Its restriction to $Z$ is the
negative of the $ij$ coordinate and $h_{ij}(c)=0$.
\end{lemma}

\begin{proof}
The class $x_i\cup x_j\in H^2(\Q,\F_3)$ is locally trivial.
At $q_i$ the unit $q_j$ is a cube because $q_i\equiv2\bmod3$,
and the same holds with $i,j$ interchanged. At $3$ both are cubes
in $\Q_3$ since $q_i,q_j\equiv-1\bmod9$. Away from these places
the two classes are unramified and the unramified Galois group has
$3$-cohomological dimension one. The real place contributes zero.
Restrict to $K$, identify $\F_3$ with $\mu_3$ there and apply the
Brauer local--global injection. Corestriction back to $\Q$ multiplies
by two, so the global cup class vanishes.
Choose continuous $1$-cochains $h_{ij}$ with $dh_{ij}=f_i\smile f_j$.
Writing $\varepsilon(g)\in\{1,-1\}$ for the cyclotomic character,
the map $g\mapsto ((f_i(g))_i,(-h_{ij}(g))_{ij},\varepsilon(g))$
is a homomorphism to $P\rtimes\langle c\rangle$ by the definition
of the inhomogeneous coboundary. Moreover $h_{ij}(c)=0$ by $c^2=1$.
The images of the $q_i$ are independent in $K^\times/K^{\times3}$
by valuations, so the image over $K$ surjects onto $W$. Commutators of
lifts of the coordinate vectors span $Z$, proving surjectivity.

We can remove auxiliary finite ramification. At a prime
$\ell\nmid3D$ the inertia image lies in $Z$, which is central even
over $\Q$. It is tame. Conjugation by Frobenius on tame inertia is
raising to the $\ell$th power, so nontrivial image forces
$\ell\equiv1\bmod3$. A cubic Dirichlet character of conductor $\ell$
generates the dual of the order-three tame inertia quotient there.
Adding suitable multiples of it to the individual $h_{ij}$ kills
the entire inertia image without changing any coboundary or introducing
ramification elsewhere. The initial homomorphism has finite image and
hence only finitely many auxiliary ramified primes. The field cut out
after these changes has the asserted properties.
\end{proof}

For the chosen \(D\) and the fixed real radicals defining the \(f_i\),
fix one such system of cochains and its field \(H\) throughout the
family. Put $H_0=K(\sqrt[3]{2},\sqrt[3]{3},\sqrt[3]{q_1},\ldots,
\sqrt[3]{q_{2m}})$ and $\mathcal H=HH_0$. The radicals occurring in $H_0$
are independent over $K$ (the valuation of $3$ over $3$ is two).
As $H_0/K$ is abelian and $[P,P]=Z$, the intersection $H\cap H_0$
is $K(\sqrt[3]{q_1},\ldots,\sqrt[3]{q_{2m}})$. Thus
$\Gal(\mathcal H/K)\simeq P\times\F_3^2$. The commutator subgroup
of $\Gal(\mathcal H/\Q)$ contains $Z$ and projects onto the whole
$W\times\F_3^2$, since $c$ acts there by negation. Consequently
$\mathcal H^{\rm ab}=K$. For $m=0$ take $H=K$ and omit the
$h_{ij}$ and $Z$ throughout.

\subsection{The comparison with the reference twist}\label{sec:variation}

For each positive $n$, write $\Theta_n$ for the Cassels--Tate pairing
on $U_n=\Sel_\phi(E_n/\Q)\subseteq H^1(\Q,\mu_3)$ with values in
$\tfrac13\Z/\Z$. In the notation of Koymans and Smith, $d=16$ gives
$E_{d,n}=A_n$, $E_{-27d,n}=E_n$, $M_d=\F_3$, and $M_{-3d}=\mu_3$.
With compatible identifications of the coefficient lines, $\Theta_n$
is their $\mathrm{CTP}_{16,n}$ on $U_n$
\cite[Sections~4.1--4.3]{KS}. For $z\in Z$ put
$b_{ij}(z)=-h_{ij}(z)$, the upper-right coordinate of the corresponding
three-by-three unipotent matrix.

\begin{lemma}\label{lem:variation}
Let \(k=Dt\), where \(t=p_1\cdots p_6\) is squarefree, prime to
\(6D\), and \(t\equiv1\bmod9\). Suppose that \(p_2,p_3\) split in
\(\Q(\zeta)\), and the four other primes are inert there. Prescribe
the rational Frobenius conjugacy classes in \(H\) to be those of
\(z,1\), with \(z\in Z\), at \(p_2,p_3\), respectively, and that of \(c\) at each inert
prime of \(t\). Then \([q_i]\in U_D\cap U_k\) for every \(i\).
There is a single \(\epsilon\in\F_3^\times\), independent of \(i,j\)
and \(k\), for which
\[
 \Theta_k([q_i],[q_j])-\Theta_D([q_i],[q_j])
    =\frac{\epsilon\,b_{ij}(z)}{3}\pmod{\Z}\qquad(i<j).
\]
An element of \(\F_3\) on the right is represented by any integer.
\end{lemma}

\begin{proof}
Write \(S_{\mathrm{fix}}=\{\infty,2,3,q_1,\ldots,q_{2m}\}\),
\(R=\{p_1,\ldots,p_6\}\), and \(S=S_{\mathrm{fix}}\cup R\).
At \(3\), each \(q_i\) is a local cube. At \(q_j\), the cube map is
bijective on units, and the local descent line generated by \(D\)
or \(k\) contains every \([q_i]\), since both \(D\) and \(k\) have
valuation one there. At \(2\) the \(q_i\) are unit cubes, and at the
real place all elements are cubes. At every \(r\in R\) all \(q_i\)
are cubes, by the Frobenius hypotheses in the split case and the
cube map on units in the inert case. At finite places outside \(S\)
these classes satisfy the unramified Kummer condition of good
reduction. This proves the assertion about \(U_D\cap U_k\).

Put \(a=\chi_{3,t}\) for the Kummer cocycle of the positive real
cube root of \(t\). At every \(v\in S_{\mathrm{fix}}\), \(t\) is a cube in
\(\Q_v\). At \(3\) this uses \(t\equiv1\bmod 9\), and at the other
finite places it follows because \(t\) is a unit and the residue
characteristic is \(2\bmod3\). For these places, choose the local
embeddings of \(\overline{\Q}\) so that the chosen cube root lies
in \(\Q_v\). Thus \(a|_{G_v}=0\) as a cocycle.
For \(r\in R\) inert, choose the local embedding whose Frobenius in
\(H\) is \(c\) itself. At the split primes the prescribed central
elements are singletons. Since \(H\) is unramified at all \(r\in R\),
we obtain \(f_i|_{G_r}=0\) as cocycles for every \(i\). These choices
of embeddings can be made independently at distinct places.

\cite[Proposition~4.1]{KS} compares \(b_0p^e\) with \(b_0q^e\) for a fixed twisting parameter \(b_0\). Here
the local identifications on \(S_{\mathrm{fix}}\) and the simultaneous
vanishing of every \(f_i\) at the six places in \(R\) allow us to
apply its cochain argument directly to \(D\) and \(Dt\). Fix \(f=f_i\) and \(g=f_j\), and use the comparison
\(b=\beta_{16,D}^{-1}\beta_{16,k}\) of
\cite[(4.1)--(4.4)]{KS}, with the cube root chosen compatibly with
that of \(t\). It identifies the two exact sequences of abelian groups
\[
 0\longrightarrow \F_3\longrightarrow A_n[3]
   \longrightarrow\mu_3\longrightarrow0,\qquad n=k,D,
\]
inducing the identity on the submodule and the quotient.
At places of \(S_{\mathrm{fix}}\), \(b\) is induced by an isomorphism of the local
curves and identifies their full local Kummer conditions. Select a
local lift \(\Phi_v\) of \(f\) satisfying this condition on \(A_k[3]\)
and use \(b\Phi_v\) for \(A_D[3]\). The representative of a local
lift can be adjusted by a coboundary to project to the chosen cocycle
\(f|_{G_v}\). At all places of \(R\) take the zero local full lift on
both curves.

Let $\sigma_0$ be the additive section for $A_k[3]$ of \cite[(4.2)]{KS},
and let $\iota(g)$ be the cocycle with values in
$\operatorname{Hom}(\F_3,\mu_3)$ induced by the Weil pairing as
in the proof of \cite[Proposition~4.1]{KS}. Let \(\Q_S\) be the
maximal extension of \(\Q\) in \(\overline{\Q}\) unramified outside
\(S\), let \(G_{\Q,S}=\operatorname{Gal}(\Q_S/\Q)\), and let
\(\mathcal O_{\Q_S}\) denote the integral closure of \(\Z\) in
\(\Q_S\). Since \(3\in S\), the field \(\Q_S\) contains
\(\Q(\zeta)\), so \(\mu_3\subseteq\mathcal O_{\Q_S}^{\times}\).
Using additive notation for cochains, take a two-cochain
\(e\in C^2(G_{\Q,S},\mathcal O_{\Q_S}^{\times})\) as in
\cite[Definition~2.2]{KS} with
\[
 de=d(\sigma_0\circ f)\smile\iota(g).
\]
Such a primitive exists by the cohomological vanishing used in
\cite[Definition~2.2 and proof of Proposition~4.1]{KS}.
Via $L\otimes L\simeq\F_3$, the coefficient isomorphism
of \cite[(4.8)]{KS} is a fixed nonzero scalar $\lambda\in\F_3^\times$
times the natural map $L\otimes\F_3\to L$. Formula~(4.3) of
\cite{KS}, applied to $k,D$, gives the cochain identity
\[
 d(b\,\sigma_0\circ f)\smile\iota(g)-d(\sigma_0\circ f)\smile\iota(g)
       =\lambda\, a\smile f\smile g .
\]
Indeed the comparison cocycle is $\chi_{3,D^2/k^2}=\chi_{3,t^{-2}}=a$
for our compatible choice of cube roots. Define
\[
 e_D=e-\lambda\,a\smile h_{ij}.
\]
Both terms are cochains in the same complex because $a$ and $h_{ij}$
factor through $G_{\Q,S}$ and their cup product takes values in
$\mu_3\subseteq\mathcal O_{\Q_S}^{\times}$. Since $da=0$ and
$dh_{ij}=f\smile g$, the cup-product differential gives
$d(-a\smile h_{ij})=a\smile f\smile g$. Thus
$de_D=d(b\,\sigma_0\circ f)\smile\iota(g)$, as required for the
reference twist.

In the Weil-pairing formula for $\Theta_k-\Theta_D$, the local-lift
terms at $S_{\mathrm{fix}}$ agree since $b$ preserves the Weil pairing.
They vanish at $R$ by our choices of local lifts and cocycles.
The signs in the formula for each pairing are those of
\cite[proof of Proposition~4.1]{KS}, so their difference is
\[
 \sum_{v\in S}\operatorname{inv}_v\bigl((e_D-e)|_{G_v}\bigr)
 =\sum_{v\in S}\operatorname{inv}_v
      \bigl((-\lambda\,a\smile h_{ij})|_{G_v}\bigr).
\]
These restrictions are two-cocycles. At each place in $S$ at least
two of the classes of $a,f,g$ vanish locally. At $q_i$ only $f$
may survive, at $q_j$ only $g$ may survive, and at the moving primes
both $f$ and $g$ vanish. At $2$, $3$, the real place, and the other
fixed primes all three vanish. Consequently every permuted pairwise
cup product is trivial locally. This verifies both local conditions of
\cite[Definition~2.2]{KS} for the generalized R\'edei symbol. The
chosen embeddings make at least one cocycle identically zero at each
such place, so the local correction cochains in that definition can
all be chosen zero.

At \(v\in S_{\mathrm{fix}}\) the cocycle \(a\) vanishes.
For \(p\in R\), the restrictions of \(f_i,f_j\) vanish as cocycles,
and \(h_{ij}|_{G_p}\) is an unramified character.
Let \(F_p\) be the arithmetic Frobenius determined by the chosen
embedding. Local reciprocity gives
\begin{equation}\label{eq:moving-invariant}
 \operatorname{inv}_p\!\left(
     -\lambda\,\operatorname{res}_{G_p}(a\smile h_{ij})\right)
   =\frac{\epsilon}{3}\,v_p(t)\bigl(-h_{ij}(F_p)\bigr)
          \quad\text{in }\tfrac13\Z/\Z .
\end{equation}
Indeed, an unramified character is trivial on units and its value on
the reciprocity image of \(t\) is \(v_p(t)\,h_{ij}(F_p)\) under the
arithmetic Frobenius convention. The coefficient map \(\lambda\)
and the ordering of the cup factors fix one
\(\epsilon\in\F_3^\times\). This calculation takes place over
\(\Q_p\), so the same scalar applies at split and inert primes.
Here \(v_p(t)=1\) for all six primes. The values of \(h_{ij}\) are
zero at the four inert primes since \(h_{ij}(c)=0\), and zero at
\(p_3\) since \(F_{p_3}=1\). At \(p_2\) one has
\(-h_{ij}(F_{p_2})=b_{ij}(z)\), which gives the claimed difference.
\end{proof}

With the cocycles \(f_i\) fixed, consider another system
\(h'_{ij}\) as in Lemma~\ref{lem:group}, inflated along with
\(h_{ij}\) to the absolute Galois group. The difference
\(h'_{ij}-h_{ij}\) has zero coboundary and hence is a
character. Both systems factor through extensions unramified outside
\(3D\), so it too is unramified outside $3D$. No \(q_i\) admits a ramified cubic character
of \(\Q_{q_i}^{\times}\), since \(q_i\equiv2\pmod3\). By
Kronecker--Weber every such global character is a multiple of
\(\chi_9\). Indeed, for \(n\ge2\) the cubic characters modulo
\(3^n\) all factor through modulus \(9\). Individual Frobenius
coordinates may therefore change. By \eqref{eq:moving-invariant}
the change in their weighted sum is a multiple of
\[
 \sum_{p\mid t}v_p(t)\chi_9(p)=\chi_9(t)=0 ,
\]
where the last equality follows from \(t\equiv1\pmod9\).
Equivalently, the change in the sum of local invariants is the
reciprocity sum of the global Brauer class
\(-\lambda\,a\smile(h'_{ij}-h_{ij})\).

For \(m=1\), the rational kernel point \(-T=(0,-4D)\) on \(A_D\) has the
radical class \([D]=[q_1]+[q_2]\), and hence
\(\Theta_D([q_1],[q_2])=0\).
For \(m\ge2\), the same argument puts the all-ones vector in the
radical of the base pairing on these fixed classes. Its other entries
remain part of the base matrix in Lemma~\ref{lem:variation}.

\begin{corollary}\label{cor:all-frob}
For every parameter satisfying \eqref{eq:fam-primes} and
\eqref{eq:fam-chars}, attach to each $p\mid t$ a central element
$z_p\in Z$ as follows. If $p$ splits in $K$, its rational Frobenius
in $H$ equals $z_p$. If $p$ is inert in $K$, its rational Frobenius
in $H$ is conjugate to $z_pc$. Then
\[
 \Theta_k(x_i,x_j)=\Theta_D(x_i,x_j)
     +\frac{\epsilon}{3}\sum_{p\mid t} b_{ij}(z_p)
          \quad\text{in }\tfrac13\Z/\Z .
\]
The sum uses integer representatives of the coordinates.
\end{corollary}

\begin{proof}
At a split prime the Frobenius in $H$ is central because the conditions
on $q_i$ make its $W$-coordinate zero. For inert $p$, let $F_p$ be any
representative. Conjugating $F_p$ by an element projecting to $w\in W$
translates its $W$-coordinate by $2w$, so a choice of $w$ makes that
coordinate zero. The new representative is $z_pc$. Its square $z_p^2$
is central and equals $F_p^2$, which determines $z_p$ uniquely.
Choose the embeddings giving these representatives. Then all the
$f_i$ vanish as local cocycles at the six moving primes, and at each of
them $h_{ij}(F_p)=h_{ij}(z_p)=-b_{ij}(z_p)$ in the chosen embedding.
Here $h_{ij}(c)=0$ and the cup correction in the product $z_pc$ is zero.
Apply \eqref{eq:moving-invariant} at each of these rational primes.
Since \(v_p(t)=1\) for every \(p\mid t\), summing gives the formula.
\end{proof}

\section{Realizing the pairing and counting twists}\label{sec:counting}
\subsection{Proof of the realization theorem}\label{sec:count}
Identify $\tfrac13\Z/\Z$ with $\F_3$ by $j/3\mapsto j$.
\begin{proof}[Proof of Theorem~\ref{thm:realization}]
For the parameters of Section~\ref{sec:family}, Propositions~\ref{prop:points}
and~\ref{prop:first} give the two independent points, the vanishing of the
dual isogeny Selmer group and the asserted basis of $U_k/W_k$.
Let $B_0=(\Theta_D(x_i,x_j))_{ij}$. Choose $z\in Z$ in
Lemma~\ref{lem:group} with coordinates
$\epsilon b_{ij}(z)=B_{ij}-(B_0)_{ij}$ for all $i<j$.
Use the coordinates of $\Gal(\mathcal H/K)=P\times\F_3^2$ from
Section~\ref{sec:group}, the last two being the cubic Kummer
coordinates of $2$ and $3$, in that order. Prescribe in
$\Gal(\mathcal H/\Q)$ the conjugacy classes of the following elements.
\[
 \begin{array}{c|ccc}
  &p_2 & p_3 & p_1,p_4,p_5,p_6 \\ \hline
  \text{representative}&(z,1,0)&(1_P,0,1)&c .
 \end{array}
\]
Here $c$ is complex conjugation with respect to the real choices of
the radicals. The projections to $H$ at the split primes are
singletons because $Z$ is central even over $\Q$. At a split prime,
choosing a prime above it giving the indicated representative
normalizes the cubic character with the values in
\eqref{eq:fam-chars}. The projections to $H$ then give $B$ by
Lemma~\ref{lem:variation}. For $m=0$ the matrix condition is empty.

We check the hypotheses of the simultaneous prime-values theorem
\cite[Theorem~1.4 and Corollary~1.6]{FT}. Set $M=90D$. Restrict
$a,b$ to residue classes $a\equiv5$, $b\equiv4\pmod9$ and
$a\equiv b\equiv1\pmod{10D}$. Write $a=a_0+Mx$, $b=b_0+My$.
The linear parts of the resulting six affine-linear forms in $x,y$
are pairwise nonproportional. At every prime dividing $M$ all six
values are units. At a prime $\ell\nmid M$ the change from $x,y$ to
$a,b$ is invertible modulo $\ell$, and the forbidden nonzero
directions lie on at most six lines. Since $\ell\ge7$ and
$\#\mathbf P^1(\F_\ell)=\ell+1$, a direction survives.
The maximal abelian subextension of $\mathcal H/\Q$ is $K$,
of conductor three. We use modulus three in the compatibility condition. Its Artin symbols agree with the residue classes
of all six forms modulo three. This verifies the abelian compatibility
condition for the entire tuple. There is a positive open cone on
which all six homogeneous linear parts are positive.

For a quantitative bound take a closed rectangle of positive area
strictly inside that cone and inside $(-1,1)^2$. On its dilation by
$N$ all six coordinates are $\asymp_D N$ for large $N$. The weighted
formula in \cite[Theorem~1.4]{FT}, with the positivity check above,
gives a sum $\gg_D N^2$. By the unweighting argument of Section~\ref{sec:prelim},
$\gg_D N^2/(\log N)^6$ pairs give six distinct prime values.
Each $k=D\prod p_i$ has at most $6!$ representations, since an
ordering recovers $a=p_1,b=p_2$. Also $k\ll_D N^6$, so a sufficiently
small constant multiple of $X^{1/6}$ suffices for $N$.
For fixed $D$ there are only finitely many matrices $B$, allowing the
implied positive constant to depend on $D$ alone.
\end{proof}

\subsection{The nondegenerate case}
\begin{proof}[Proof of Theorem~\ref{thm:main}]
Take a block diagonal alternating matrix of size $2m$ with $m$ copies
of $\left(\begin{smallmatrix}0&1\\-1&0\end{smallmatrix}\right)$.
For $m=0$ use the unique form on the zero space. Define
$\mathcal K_{2m}(D)$ by \eqref{eq:fam-primes}, \eqref{eq:fam-chars}
and the requirement that the matrix on the marked quotient be this
chosen matrix. The construction in the proof of Theorem~\ref{thm:realization} gives its counting
bound. By Proposition~\ref{prop:first}, the three rational classes
span $W_k$ and the $[q_i]$ project to a basis of $U_k/W_k$.
Thus $d=2m=\dim_{\F_3}U_k-3$.
Theorem~\ref{thm:criterion} therefore gives exact rank two,
$G_k\simeq(\Z/3\Z)^{2m}$ and $G'_k=0$.

Finally an isomorphism of the short models over $K$ would give
$u^6=(k/k')^2$ for some $u\in K^\times$. Hence $k/k'$ has a cube
root in $K$, since a possible minus sign is a rational cube.
A rational noncube has irreducible cubic polynomial $X^3-k/k'$ and
cannot have a cube root in a quadratic field. Distinct positive
squarefree parameters have noncubic ratio.
\end{proof}

\subsection{A nonzero example with two fixed primes}\label{sec:explicit}

Take $q=17$ and $q'=233$, both prime and congruent to $-1\pmod 9$.
Let $F=\Q(\zeta)$ with $\zeta$ a primitive cube root of unity, and take
the positive real roots $\rho^3=17$, $r^3=233$. The identity $6^3+17=233$
gives $N_{F(\rho)/F}(6+\rho)=233$. Put
\[
 \xi=6+\rho,\qquad \alpha=6+\zeta \rho,\qquad \beta=6+\zeta^2\rho,
 \qquad w^3=\theta=\alpha\beta^2,\qquad H=F(\rho,r,w).
\]
The two radicands $17,233$ are independent modulo $F^{\times3}$ by
their valuations. Over $F$, the automorphisms $\sigma:\rho\mapsto\zeta \rho$
and $\tau:r\mapsto\zeta r$ of $F(\rho,r)$ each fix the other radical.
Since $\sigma\theta/\theta=\xi^3/233$, they extend to $H$ by
\[
   \sigma w=(\xi/r)w,\qquad \tau w=w .
\]
Indeed $w$ cannot belong to $F(\rho,r)$. Otherwise its transforms under
$\sigma$ and $\tau$ would have these values up to cube roots of unity,
and $\sigma\tau w=\zeta\tau\sigma w$ would contradict commutativity.
The two extensions have order three and satisfy
$[\sigma,\tau]=z$, where $z(w)=\zeta w$ and $z$ fixes $F(\rho,r)$.
Thus $\Gal(H/F)$ is the Heisenberg group of order $27$ with center
$\langle z\rangle$. Complex conjugation $c$ fixes $\rho,r$, interchanges
$\alpha,\beta$, and sends $w$ to $\alpha\beta/w$.
Hence $H/\Q$ is Galois, with $c\sigma c^{-1}=\sigma^{-1}$,
$c\tau c^{-1}=\tau^{-1}$ and $czc^{-1}=z$.
Only $3,17,233$ can ramify, since $\xi\alpha\beta=233$ and all three
factors are algebraic integers.
With the coordinates on the two radicals ordered as $17,233$, the
commutator $z$ is the central element of coordinate one in
Lemma~\ref{lem:group}. We use this choice of governing field for the
following calculation.

\begin{example}\label{ex:nonzero-fiber}
For $m=1$, take $D=17\cdot233=3961$ and $(a,b)=(5146061,17488921)$. Set
$u=213894000a$ and
\[
 k=3961\,ab(u+b)(u-b)(2u-b)(u/2-b).
\]
For this parameter,
\[
 \operatorname{rank}E_k(\Q)=2,\qquad
 \Sha(E_k/\Q)[3^\infty]\simeq(\Z/3\Z)^2,\qquad
 \Sha(A_k/\Q)[3^\infty]=0 .
\]
\end{example}

\begin{proof}
Here $u=c_0a$. In the order of \eqref{eq:fam-primes} the six moving
factors are distinct primes, and the two split primes $p_2,p_3$
satisfy \eqref{eq:fam-chars}. Appendix~\ref{app:nonzero-fiber}
gives the primality and residue computations. In the explicit field
$H$ above, write $z(w)=\zeta w$ and let $c$ denote complex
conjugation. Representatives of the six arithmetic Frobenius classes
are, in order,
\[
                   (c,z^2,1,c,c,c).
\]
On the reference twist, $[D]=[17]+[233]$ is the Kummer class of
the rational point $(0,-4D)$ on $A_D$. It lies in the radical of
$\Theta_D$. Alternation therefore gives
$\Theta_D([17],[233])=\Theta_D([17],[D])=0$.
In the central coordinate $z\mapsto1$, Lemma~\ref{lem:variation}
now gives
\[
 \Theta_k([17],[233])=\frac{2\epsilon}{3}\ne0
       \quad\text{in }\tfrac13\Z/\Z ,
\]
where $\epsilon\in\{1,-1\}$ accounts for the common pairing and
reciprocity conventions. Either choice gives a nondegenerate
alternating form on the two-dimensional quotient with basis
$[17],[233]$ from Proposition~\ref{prop:first}. Theorem~\ref{thm:criterion}
gives the asserted rank and both three-primary groups.
\end{proof}

\section{A disjoint vanishing family}\label{sec:six-split}
When the rational classes exhaust $U_k$, the first descent already
determines the whole three-primary group. We construct a second family
with four fixed inert prime factors and six moving split primes.
\subsection{The descent matrix}
\begin{definition}\label{def:pattern}
Fix
\[
 (q_1,q_2,q_3,q_4)=(17,53,71,89),\qquad
 (c_1,\ldots,c_6)=(2^4 3^3 5^3q_1,1,q_2,q_3q_4,1,2).
\]
Require $n_1>n_2>0$ and $L_i=c_ip_i$ for six distinct primes
$p_i\equiv4\pmod9$. Put $(r_1,\ldots,r_7)=(q_1,q_2,q_3,q_4,2,3,5)$.
For each $i$ require a primitive cube root $\omega_i\in\F_{p_i}$
such that the additive cubic character defined by
$x^{(p_i-1)/3}=\omega_i^{\chi_i(x)}$ satisfies
$\chi_i(r_j)=\delta_{ij}$ for $1\le j\le7$.
The associated parameter is $k=q_1q_2q_3q_4\prod_i p_i$.
\end{definition}
Write $\mathcal V$ for the set of parameters in Definition~\ref{def:pattern}.

Here $d_i=c_i$ and $2\prod_i c_i=60^3D$, so
\eqref{eq:sextic-points} holds with $\varepsilon=1$.
The residue of $k$ modulo nine is $(-1)^44^6=1$.

\begin{theorem}\label{thm:fiber}
For every parameter of Definition~\ref{def:pattern},
\[
 \dim_{\F_3}\Sel_\phi(E_k/\Q)=3,\qquad
 \Sel_{\widehat\phi}(A_k/\Q)=0,\qquad
 \dim_{\F_3}\Sel_3(E_k/\Q)=2 .
\]
In particular, $\operatorname{rank}E_k(\Q)=2$ and
$\Sha(E_k/\Q)[3^\infty]=0$.
\end{theorem}

\begin{proof}
For $i\le4$, the character $\chi_i$ kills $-1,2,3,5$ and takes
$q_j$ to $\delta_{ij}$. The absolute determinants in
Lemma~\ref{lem:eliminate} are $1,2,3$. The unique $c_j$ divisible by
$q_i$ has valuation one there and index $j=\nu(i)$, where
$(\nu(1),\nu(2),\nu(3),\nu(4))=(1,3,4,4)$.
The four transformed rows are $e_{q_i}-e_{p_{\nu(i)}}$,
with $e_v$ the unit row at $v$. If $\nu(i)=i$, the diagonal entry
comes from the row sum.

At $p_5$ only $2$ has nonzero character among the $r_j$.
For $j=1,2,3,4,6$ the determinants are $-1,2,-3,1,3$, and the
$2$-adic valuations of $c_j$ are $4,0,0,0,1$.
The resulting off-diagonal entries are $2,1,0,0,2$.
At $p_6$ only $3$ has nonzero character. For $j=1,2,3,4,5$
the determinants are $-2,1,-3,-1,-3$, and the $3$-adic valuations of
$c_j$ are $3,0,0,0,0$. This gives $0,0,1,0,1$ off the diagonal.
The diagonal entry in each of these two rows is $1$. Together with $w$,
they give the three rows of $C$. Thus
\[
 M\ \sim\
 \begin{pmatrix}I_4&R\\0&C\end{pmatrix},\qquad
 R=\begin{pmatrix}
 2&0&0&0&0&0\\
 0&0&2&0&0&0\\
 0&0&0&2&0&0\\
 0&0&0&2&0&0
 \end{pmatrix}.
\]
The columns $q_1,\ldots,q_4,p_2,p_3,p_4$ give a minor of determinant
$1$, so $M$ has rank seven. The fixed-prime parts of the three
rows in \eqref{eq:point-vectors} are respectively
$(1,1,1,1)$, $(1,1,2,2)$ and $(0,2,1,1)$, by the expressions for
$\delta(P)/[k]$ and $\delta(Q)/[k]$. These three independent rational
classes form a basis of the right kernel, so $W_k=\Sel_\phi(E_k/\Q)$.
The left kernel is zero. Proposition~\ref{prop:descent} gives the three
Selmer dimensions. Theorem~\ref{thm:criterion} applies with $a=3$ and $d=0$
and gives rank two and the vanishing of both primary groups.
\end{proof}

\subsection{Counting the ten-prime twists}
\begin{proposition}\label{prop:prime-count}
For all sufficiently large $X$, there are $\gg X^{1/3}/(\log X)^6$
distinct integers $k\le X$ satisfying the prime pattern of
Definition~\ref{def:pattern}.
\end{proposition}

\begin{proof}
We first impose congruences on \(m,n\). In the last column of the table
we list all positive valuations among the six forms at the prime
underlying the modulus.
\[
\begin{array}{c|c|l}
 \text{modulus}&(n_1,n_2)\ \text{modulo that modulus}&\text{positive valuations}\\
 \hline
 2^5 &(16,1)&v_2(L_1)=4,\quad v_2(L_6)=1\\
 3^5 &(4c_1,4)&v_3(L_1)=3\\
 5^4 &(125,1)&v_5(L_1)=3\\
 q_1^2 &(q_1,1)&v_{q_1}(L_1)=1\\
 q_2^2 &(1+q_2,1)&v_{q_2}(L_3)=1\\
 q_j^2,\ j=3,4 &(q_j-1,1)&v_{q_j}(L_4)=1
\end{array}
\]
The entries follow by substitution. The Chinese remainder theorem
gives a solution \((n_1^*,n_2^*)\) modulo
\[
 W=2^5 3^5 5^4\prod_{j=1}^4 q_j^2.
\]
Write \((n_1,n_2)=(n_1^*+Wx,n_2^*+Wy)\). The quotients
\(\psi_i(x,y)=L_i(n_1,n_2)/c_i\) are integral affine-linear forms since
\(c_i\mid W\) and \(c_i\mid L_i(n_1^*,n_2^*)\). The table makes every
\(\psi_i\) a unit at every prime dividing \(W\). At \(3\), the
congruence on \(n_1\) gives \(n_1/c_1\equiv4\bmod 9\) after cancellation
of \(3^3\). For \(i=2,\dots,6\), the residues of \(L_i\) modulo \(9\)
are respectively \(4,-4,4,4,8\), and those of \(c_i\) are
\(1,-1,1,1,2\). Thus all six \(\psi_i\) are \(4\bmod 9\).

Let \(K=\mathbb Q(\zeta_3)\), and in a fixed embedding in \(\mathbb C\)
take the real cube roots in
\[
 H=K(\sqrt[3]{r_1},\ldots,\sqrt[3]{r_7}).
\]
The classes of the \(r_j\) are independent in \(K^\times/K^{\times3}\).
Indeed, the valuation of \(r_j\) at a prime above it is \(1\) or \(2\),
and all the other \(r_h\) have valuation zero there. By Kummer theory,
\(\operatorname{Gal}(H/K)=\mathbb F_3^7\) with the coordinate \(u_j\)
of \(u\) determined by
\(u(\sqrt[3]{r_j})=\zeta_3^{u_j}\sqrt[3]{r_j}\).
The field \(H\) is Galois over \(\mathbb Q\), and complex conjugation
acts on \(\mathbb F_3^7\) by negation. Its commutator with \(u\)
is \(-2u\), so the commutator subgroup of \(\operatorname{Gal}(H/\mathbb Q)\)
is precisely \(\operatorname{Gal}(H/K)\). The maximal abelian subfield
is therefore \(K\), of conductor \(3\).

For \(1\le i\le6\), let \(C_i=\{e_i,-e_i\}\), where \(e_i\) is the
\(i\)-th standard basis vector of \(\mathbb F_3^7\). These are distinct
conjugacy classes in \(\operatorname{Gal}(H/\mathbb Q)\). If an
unramified prime \(p\) has Frobenius class \(C_i\), it splits in \(K\).
Choose a prime of \(H\) above \(p\) with arithmetic Frobenius \(e_i\),
and let \(\omega\in\mathbb F_p\) denote the reduction of \(\zeta_3\).
Reduction of the Frobenius action on each cube root gives
\[
 r_j^{(p-1)/3}=\omega^{\delta_{ij}}\quad\text{in }\mathbb F_p.
\]
Consequently the character defined by \(a^{(p-1)/3}=\omega^{\chi_i(a)}\)
has all the required values.

We check the hypotheses giving a positive main term in
Frei and Tinkov\'a's prime-values theorem
\cite[Theorem~1.4 and Corollary~1.6]{FT}.
The linear parts of the \(\psi_i\) are nonzero multiples of the six
pairwise nonproportional coefficient vectors of the \(L_i\). There is
no fixed prime divisor of \(\prod_i\psi_i\). For primes dividing \(W\)
this follows from the table. For a prime \(\ell\nmid W\), the change
from \((x,y)\) to \((n_1,n_2)\) is invertible modulo \(\ell\), all
the \(c_i\) are units, and the six equations \(L_i=0\) remove at most
six directions in \(\mathbb P^1(\mathbb F_\ell)\). Since
\(\ell\ge7\), at least one direction remains. The abelian compatibility
condition holds with modulus \(3\) because every \(C_i\) fixes \(K\)
and every \(\psi_i\) is \(1\bmod3\). Finally, on the open cone
\(y>0,\ 2y<x<3y\), the six linear parts are positive.

Choose a fixed closed rectangle of positive area inside this cone and
inside $(-1,1)^2$. On its dilation by $Y$, the weighted theorem gives
a count $\gg Y^2$. The unweighting argument in
Section~\ref{sec:prelim} gives $\gg Y^2/(\log Y)^6$ pairs. The six
classes are disjoint; in particular they label the six prime factors.
These primes differ from the $q_j$ by their residues modulo nine.

For each pair, the resulting \(k\) is positive, squarefree, and
\(k\ll Y^6\). Its six moving primes are uniquely labelled by their
Frobenius classes in \(H\), and \(p_1,p_2\) determine
\(n_1=c_1p_1\) and \(n_2=p_2\). Thus distinct pairs produce distinct
\(k\). Taking \(Y\) to be a sufficiently small fixed multiple of
\(X^{1/6}\) proves the proposition.
\end{proof}

Each $k\in\mathcal V$ has six split prime divisors, while every
$k\in\mathcal K_{2m}(D)$ has two. The families are disjoint and give
nonisomorphic curves even over $K$, by the squarefree-ratio argument
in the proof of Theorem~\ref{thm:main}.

\section{The isogeny partner and the CM field}\label{sec:consequences}
Put $K=\Q(\zeta_3)$, $\mathcal O_K=\Z[\zeta_3]$ and
$\mathfrak p=(1-\zeta_3)$. The following corollary applies to the
constructions of Theorem~\ref{thm:main} and Section~\ref{sec:six-split}.

\begin{corollary}\label{cor:cm-field}
On $\mathcal K_{2m}(D)$, and on $\mathcal V$ with $m=0$, one has
\[
\begin{aligned}
 \operatorname{rank} A_k(\Q)&=2, & \Sha(A_k/\Q)[3^\infty]&=0,\\
 \dim_{\F_3}\Sel_3(E_k/\Q)&=2m+2, &
 \dim_{\F_3}\Sel_3(A_k/\Q)&=3.
\end{aligned}
\]
Moreover
\[
 \operatorname{rank}_{\Z} E_k(K)=4,\qquad
 \Sha(E_k/K)[3^\infty]\simeq(\mathcal O_K/\mathfrak p)^{2m}
 \quad\text{as $\mathcal O_K$-modules}.
\]
Restriction identifies the underlying group with $G_k$, and complex
conjugation fixes it elementwise.
\end{corollary}
\begin{proof}
The partner has the same rational rank by isogeny. In both constructions
$V_k=0$ and the isogeny pairing has maximal rank on $U_k/W_k$.
Thus Theorem~\ref{thm:criterion} gives $G'_k=0$, and
\eqref{eq:dimensions} gives the Selmer dimensions, using
$E_k(\Q)[3]=0$ and $A_k(\Q)[3]=\langle T\rangle$.

The quadratic twist $E_k^{(-3)}$ has constant coefficient
$3^6\cdot16k^2$ and is isomorphic over $\Q$ to $A_k$ by
$(x,y)\mapsto(x/9,y/27)$. The eigenspaces of complex conjugation on
$E_k(K)\otimes\Q$ therefore each have dimension two. On the
three-primary subgroup the idempotents $(1\pm c)/2$ for complex
conjugation $c$ define complementary direct summands. Restriction
and corestriction identify the invariant subgroup with the group
over $\Q$, since the extension has degree two. The same identities
locally, at places of degree one or two, ensure that everywhere
locally trivial classes correspond. Twisting applies this argument
to the minus eigenspace as well. Thus
\[
 \Sha(E_k/K)[3^\infty]^+\simeq G_k,\qquad
 \Sha(E_k/K)[3^\infty]^-\simeq G'_k=0 .
\]
For any class $x$ in this group, both $x$ and $\zeta_3x$ are fixed
by $c$, while $c(\zeta_3x)=\zeta_3^2c(x)$. Consequently
$(\zeta_3-\zeta_3^2)x=0$. Since $\zeta_3-\zeta_3^2$ generates
$\mathfrak p$ and $\mathcal O_K/\mathfrak p=\F_3$, the claimed
CM module structure follows.
\end{proof}

For the parameters in both families the discriminant of the
short model has valuation nine at $3$. The ramified quadratic base
change to $K$ doubles it to eighteen. A Weierstrass change of variables
alters it by a multiple of twelve, so the unique prime of $K$ above
three still has bad reduction. The corollary requires only
isogeny descent and restriction--corestriction at this prime.

\appendix
\section{Arithmetic for the nonzero example}\label{app:nonzero-fiber}

For Example~\ref{ex:nonzero-fiber}, the table gives the value of each
linear form, the bound used in trial division, and a choice of the
residue $S$ of $\rho=\sqrt[3]{17}$. Each entry satisfies $S^3=17$
in $\F_{p_i}$.
\[
\begin{array}{c|rcrr}
i & p_i & p_i\bmod9 & \lfloor\sqrt{p_i}\rfloor & S \\ \hline
1 & 5146061          & 5 & 2268     & 4987116          \\
2 & 17488921         & 4 & 4181     & 16138478         \\
3 & 1100711589022921 & 4 & 33176973 & 946664950548556  \\
4 & 1100711554045079 & 5 & 33176973 & 24157093138025   \\
5 & 2201423125579079 & 5 & 46919325 & 1634432193577240 \\
6 & 550355768278079  & 5 & 23459662 & 428915292745235
\end{array}
\]
Sieving through $46919325$ gives the primes for trial division.
For the six rows in order, the numbers of primes up to the indicated
bounds are $336$, $574$, $2\,041\,862$, $2\,041\,862$, $2\,827\,126$,
and $1\,475\,340$.
Dividing $p_i$ by every prime up to its bound gives a nonzero remainder
in every case. The six primes are distinct and prime to $30D$. In
particular all are unramified in $H$, whose possible ramified primes
are $3,17,233$.

For the two split primes put $e=(p_i-1)/3$ and use the indicated
reduction $Z$ of $\zeta$, so $Z^2+Z+1=0$ in $\F_{p_i}$. Both
$17^e$ and $233^e$ equal $1$ in each row. Thus each prime splits
completely in $F(\rho,r)$, where $F=\Q(\zeta)$ and $r^3=233$.
With $\theta=(6+ZS)(6+Z^2S)^2$ the remaining powers are
\[
\begin{array}{c|rr|ccc}
i& Z&\theta & 2^e & 3^e &\theta^e\\ \hline
2& 9795391         &4458251         &Z&1&Z^2\\
3& 566438750108624 &155948641933023 &1&Z&1
\end{array}
\]
Writing $x^e=Z^{\chi_i(x)}$ gives the characters in
\eqref{eq:fam-chars}. Arithmetic Frobenius acts on $w$ by the ratio
$w^{p_i}/w=\theta^e$ in the residue field. The two Frobenius elements in
$\Gal(H/F(\rho,r))=\langle z\rangle$ are consequently $z^2$ and $1$,
respectively.

For each inert prime $p=p_i$, take $S=17^{(2p-1)/3}\bmod p$ as in the
first table. Cubing permutes $\F_p$, so $r$ also reduces to an element
of $\F_p$. The residue field for $F(\rho,r)$ at these choices is
$\F_{p^2}=\F_p[Z]/(Z^2+Z+1)$. There one has
\[
 \theta=(36-6S+S^2)(6-S-SZ).
\]
For direct exponentiation, the product of $A+BZ$ and $C+EZ$ has
coordinates $(AC-BE,AE+BC-BE)\bmod p$. Repeated squaring gives
$\theta^{(p^2-1)/3}=(1,0)$ for $i=1,4,5,6$. The square of the
rational Frobenius fixes $F(\rho,r)$ and acts on $w$ by this factor,
so it is the identity. Frobenius acts nontrivially on $F$ and
therefore has order two. Since $|\Gal(H/\Q)|=54$, its subgroup of
order two is conjugate to $\langle c\rangle$ by the Sylow theorem.

The same split-prime powers give the descent matrix directly.
For $i=2,3$ put $p=p_i$, evaluate $p_j^{(p-1)/3}$ for $j\ne i$ and choose the
diagonal entry so that each row sums to zero. The two fixed columns
vanish, and the moving parts of these rows are
\[
 (0,2,1,1,2,0),\qquad (1,1,2,1,2,2).
\]
Subtracting the row $(1,1,1,1,1,1)$ at $3$ from each gives the
first two rows of $C$ in \eqref{eq:moving-block}.
Thus $\dim_{\F_3}U_k=5$, $V_k=0$ and
$\dim_{\F_3}\Sel_3(E_k/\Q)=4$ in this example.

\section{A vanishing example with nonzero dual Selmer group}\label{app:1729}
\begin{example}\label{ex:1729}
The general rank bound also applies when $V_k\ne0$. For $k=1729=7\cdot13\cdot19$, the two representations
$1^3+12^3=9^3+10^3=k$ give $P=(1,12)$ and $Q=(9,10)$ on the affine
Fermat model. Taking discrete logarithms modulo $3$ with primitive
roots $3,2,2$ at the three primes, in that order, gives
\[
 M=\begin{pmatrix}1&0&2\\2&2&2\\0&2&1\\2&1&0\end{pmatrix}.
\]
Here $M$ has rank two. The two isogeny Selmer dimensions are therefore
one and two, and \eqref{eq:sel-sequence} bounds the full Selmer
dimension by $1+2-1=2$.

For completeness, the two points give independent classes modulo $3$.
At the good primes $\ell=5,11$, cubing permutes $\F_\ell$ and
$\#E_k(\F_\ell)=\ell+1$. Reduction followed by multiplication by
$(\ell+1)/3$ thus factors through $E_k(\Q)/3E_k(\Q)$. With the
coordinates of Section~\ref{sec:prelim}, the images are
\[
\begin{array}{c|cc}
 \ell & P & Q\\ \hline
 5  &(2,4)&(2,4)\\
 11 & O   &(5,6)
\end{array}
\]
A relation $aP+bQ=0$ modulo $3E_k(\Q)$ consequently requires first
$b=0$ and then $a=0$ modulo $3$. Since $E_k(\Q)[3]=0$,
\eqref{eq:dimensions} gives rank two and $\Sha(E_k/\Q)[3^\infty]=0$.
\end{example}

\end{document}